\documentclass[journal,twoside,web]{ieeecolor}
\usepackage{generic}
\usepackage{cite}
\usepackage{amsmath,amssymb,amsfonts}
\usepackage{algorithmic}
\usepackage{graphicx}
\usepackage{algorithm,algorithmic}
\usepackage{hyperref}
\hypersetup{hidelinks=true}
\usepackage{textcomp}
\usepackage{amsmath} 
\let\labelindent\relax
\usepackage{blindtext,graphicx,psfrag,color} 
  \graphicspath{ {./images/} }
\usepackage{amsmath}               
  \allowdisplaybreaks[1]           
\usepackage{amssymb}               
\usepackage{url}                   
\usepackage{rotating}              
\usepackage{multirow}             
\usepackage{lscape}  
\usepackage{ragged2e}
\usepackage{cuted}
\usepackage{tabularx}              
\usepackage{verbatim}              
\usepackage{footnote}               
\usepackage{float}                 
\usepackage{stfloats}
\usepackage{booktabs}              
\usepackage{lipsum}                
\usepackage{enumitem}
\usepackage{subcaption}      
\usepackage{multicol}

\newtheorem{thm}{Theorem}
\newtheorem{lem}{Lemma}

\newtheorem{cor}{Corollary}

\newtheorem{rem}{Remark}
\newtheorem{assu}{Assumption}

\newcommand{\eigLambda}{\boldsymbol{\lambda}}

\def\BibTeX{{\rm B\kern-.05em{\sc i\kern-.025em b}\kern-.08em
    T\kern-.1667em\lower.7ex\hbox{E}\kern-.125emX}}
\begin{document}
\title{Solvability of the Inverse Optimal Control problem based on the minimum principle}
\author{Afreen Islam, Guido Herrmann, \IEEEmembership{Senior Member, IEEE},  and Joaquin Carrasco, \IEEEmembership{Member, IEEE}
\thanks{The work was supported by the PhD Scholarship from the Electrical and Electronic Engineering Department, the University of Manchester awarded to Afreen Islam. \textit{(Corresponding author: Afreen Islam.)}}
\thanks{A. Islam, G. Herrmann and J. Carrasco are with the Control Systems and Robotics Group, Department of Electrical and Electronic Engineering, The University of Manchester (email: afreen.islam@manchester.ac.uk, guido.herrmann@manchester.ac.uk, joaquin.carrasco@manchester.ac.uk).}} 

\maketitle

\begin{abstract}
In this paper, the solvability of the Inverse Optimal Control (IOC) problem based on two existing minimum principal methods, is analysed. The aim of this work is to answer the question regarding what kinds of trajectories, that is depending on the initial conditions of the closed-loop system and system dynamics, of the original optimal control problem, will result in the recovery of the true weights of the reward function for both the soft and the hard-constrained methods \cite{c20}, \cite{c10a}.  Analytical conditions are provided which allow to verify if a trajectory is sufficiently conditioned, that is, holds sufficient information to recover the true weights of an optimal control problem. It was found that the open-loop system of the original optimal problem has a stronger influence on the solvability of the Inverse Optimal Control problem for the hard-constrained method as compared to the soft-constrained method. These analytical results were validated via simulation.
\end{abstract}


\section{Introduction}
Inverse Optimal Control (IOC) is an approach to learn the unknown cost function of optimal processes, such as from natural movements or locomotion in humans, animals and birds. In Optimal Control, the weights imposed on the different states and controls in a cost function are manually tweaked by trial and error methods, to generate a desired system behaviour \cite{ca}. By different permutations and combinations of penalties imposed on states and controls, weights are picked which somehow fulfill the objectives of the desired control problem. There is no systematic and scientific way to pick these penalties. However, in recent times, the research direction has shifted towards developing systematic ways to find these weights, leading to the emergence of the field of IOC. In the IOC problem, systems are observed and their dynamics are replicated under the assumption that they were generated following an optimality criterion. Hence, it is assumed that the states and control trajectories are measurable and are a result of minimizing an objective function \cite{c20}. Based on these measurements, the aim is to recover the unknown reward function of the original optimal control problem \cite{c20}. 

Kalman was the first to state the IOC problem for a linear system in \cite{c1}. Initial research was carried out predominantly by the Machine Learning community in a Markov Decision Process (MDP) setting \cite{c2,c3,c4,c5}. Mombaur proposes a bi-level approach to identify a reward function to replicate human motion in a humanoid robot \cite{c19}.  In recent years, the IOC problem has been formulated as a control systems problem. Keshavarz \cite{c7a} and Johnson \cite{c20} formulates this problem by exploiting the necessary conditions of optimality. The works in \cite{c10a,c9,c10} are also based on the minimum principle method. The works in \cite{c12,c14,c15,c18} are based on Adaptive Dynamic Programming (ADP) based approaches.

Molloy utilises the initial value of the Riccati solution associated with the IOC problem of the soft-constrained method to find conditions under which the IOC method will give unique solutions \cite{c10a}. However, to the best of knowledge of the authors, there is no work in the existing literature which mentions the kind of observations, initial conditions of the closed-loop system and the system dynamics, of the original optimal control problem, for which the IOC problem based on the minimum principle will be solvable. In this work, the solvability of the IOC problem will be analysed for different kinds of trajectories, initial conditions of the closed-loop system and the open-loop system dynamics for the original optimal control problem, for both the soft and the hard-constrained methods \cite{c20}, \cite{c10a}. The analysis will eventually be narrowed down to second-order systems because they can be seen as building blocks of higher-order systems.

The paper is organised as follows. Section II covers the preliminaries of the minimum principle based soft- and hard-constrained IOC methods. Section III deals with solvability of the IOC problem using both the methods for a general nonlinear system in the original optimal problem. Section IV considers the solvability analysis for a Linear Time Invariant (LTI) system in the original optimal problem. Section V covers examples and Section VI concludes the paper.

\section{Preliminaries} 
In this section, the preliminaries of the minimum principle based soft-constrained \cite{c7a, c20}, and hard-constrained IOC methods \cite{c10a} will be covered. The following notations will be used: $\text{Re}(\Xi)$ and $\Xi^*$ denote respectively the real part and the conjugate of a complex number/matrix $\Xi$. For a given vector $\Xi$, $\hbox{diag}(\Xi)$ is a diagonal matrix whose diagonal elements are formed by the elements of the vector $\Xi$.  $\text{ker}(\Xi)$, $\Xi^T$, $\begin{vmatrix} \Xi
\end{vmatrix}$ and $\eigLambda(\Xi)$ denote respectively the kernel, the transpose, the determinant and a set consisting of all the eigenvalues of a matrix $\Xi$. $I$ and $\mathbf{0}$ denote respectively an identity matrix and a zero vector/matrix of appropriate dimensions. 

In the spirit of \cite{c20}, the following class of optimal control problems is considered:
\begin{equation}
\begin{aligned}
& \underset{u(t)}{\text{minimize}}
& &  \int_{t_0}^{t_f} c^T\phi(t,x(t),u(t)) dt \\
& \text{subject to}
& & \dot{x}(t)=f(t,x(t),u(t))\\
& & & x(0)=x_0
\label{original_optimal_problem}
\end{aligned}
\end{equation}
where $x(t) \in \mathcal{X} \subset \mathbb{R}^{n}$ is the state, $u(t) \in \mathcal{U} \subset \mathbb{R}^{m}$ is the input, $\phi:\mathbb{R}\times \mathcal{X} \times \mathcal{U} \rightarrow \mathbb{R}^{k}_+$ are known basis functions and for the IOC problem, $c \in \mathbb{R}^{k}$ is the vector to be found, which consists of the unknown weights imposed on the states and control signals in the cost function. It is assumed that one of the elements of the vector, $c$, is known. 

\begin{rem} 
Here in this paper, the problem in (\ref{original_optimal_problem}) is referred to as the Primary Optimal Control problem.
\end{rem}

\begin{assu}
It is assumed that $f \in C^1$, that is, continuously differentiable \cite{c20}.
\label{assu_diff}
\end{assu}
Assumption \ref{assu_diff} results in the system,
\begin{equation}
\dot{x}(t) = f(t,x(t), u(t)) ,
\end{equation}
to be well-posed, that is, for every initial condition, $x_0$,
and every admissible control, $u(t)$, the system, $\dot{x}(t)=f(t,x(t), u(t))$, has a unique solution, $x(t)$, with $t \in [t_0, t_f ]$.

\begin{assu}
It is assumed that the basis function, $\phi(t,x(t),u(t))$ is known. However, the weights, $c$, of the cost function are unknown.
\end{assu}
The Hamiltonian function for the problem is:
\begin{align}
H(t,x(t),u(t),p(t))&=c^T\phi(t,x(t),u(t))\notag\\&+p^T(t)f(t,x(t),u(t)),
\label{hamiltonian}
\end{align}
where $p(t)$ is the costate trajectory. If $(x(t),u(t))$ are optimal, that is, $(x(t),u(t))=(x^*(t),u^*(t))$ where $(x^*(t),u^*(t))$ is the optimal solution of (\ref{original_optimal_problem}) for a given $c$, then there exists a costate trajectory, $p(t):\mathbb{R}\rightarrow\mathbb{R}^n$ (here, $p(t)=p^*(t)$, where $p^*(t)$ is the optimal costate trajectory) and it is necessary that
\begin{align}
\dot{p}(t)^T+\nabla_x H(t,x(t),u(t),p(t))=\mathbf{0} \label{first_cond}\\
\nabla_u H(t,x(t),u(t),p(t))=\mathbf{0}. \label{second_cond}
\end{align}
The following condition together with (\ref{first_cond}) and (\ref{second_cond}) is sufficient and necessary for optimality: $\frac{\partial^2{H(t,x(t),u(t),p(t))}}{\partial{u^2}} > \mathbf{0}$.

\subsection{Soft-constrained method}
In this subsection, the soft-constrained IOC method \cite{c20} based on the minimum principle will be described briefly.  In \cite{c20} and \cite{c7a}, the idea that the conditions in (\ref{first_cond}) and (\ref{second_cond}) will hold approximately for an approximately optimal system, was exploited to propose the soft-constrained IOC method. The left-hand sides of (\ref{first_cond}) and (\ref{second_cond}) are arranged into a residual function, $r(z(t),v(t))$, defined as:
\begin{align}
r(z(t),v(t))=&\begin{bmatrix}
\nabla_x\phi\vert_{(x,u)}^T & \nabla_xf\vert_{(x,u)}^T \\
\nabla_u\phi\vert_{(x,u)}^T & \nabla_uf\vert_{(x,u)}^T
\end{bmatrix}z(t)+ \begin{bmatrix}
          I\\
          \mathbf{0}
         \end{bmatrix}v(t)\notag\\=&F(t)z(t)+G(t)v(t),
         \label{eqn:residual}
\end{align}
where $z(t)=\begin{bmatrix}
          c\\
          p(t)
         \end{bmatrix} \in \mathbb{R}^{k+n}$ and $v(t)=\dot{p}(t)\in \mathbb{R}^n$. The following optimization problem is formulated \cite{c20} to recover the unknown weights, $c$:
\begin{equation}
\begin{aligned}
& \underset{z(t),v(t)}{\text{minimize}}
& &  \int_{t_0}^{t_f} \Vert r(z(t),v(t)) \Vert^2 dt \\
& \text{subject to}
& & \dot{z}(t)=\begin{bmatrix}
             \mathbf{0}\\
             I
\end{bmatrix}v(t) \\
& & & z(0)=z_{0} \hspace{0.1cm}\mbox{(unknown).}
\end{aligned}
\end{equation}
This is equivalent to the following Linear Quadratic Regulator (LQR) problem \cite{c20}:
\begin{equation}
\begin{aligned}
& \underset{z(t),v(t)}{\text{minimize}}
& \int_{t_0}^{t_f} &(z^T(t)Q_1(t)z(t)+v^T(t)R(t)v(t)\\
& &     &+2z^T(t)S(t)v(t)) dt \\
& \text{subject to}
& & \dot{z}(t)=A_1(t)z(t)+B(t)v(t), \\
& & & z(0)=z_{0} ,
\end{aligned}
\label{eqn:secondary_lqr}
\end{equation}
where $A_1(t)=\mathbf{0}_{(n+k)\times(n+k)}$, $B(t)=B=\begin{bmatrix}
             \mathbf{0}_{k\times n}\\
             I_{n \times n}
\end{bmatrix}$, $Q_1(t)=F(t)^TF(t)$,  $R(t)=G(t)^TG(t)$ and $S(t)=F(t)^TG(t)$. 
\begin{rem}
Here in this paper, the time varying LQR form of the residual optimization problem in (\ref{eqn:secondary_lqr}), that is, the IOC problem, is termed as the Secondary Optimal Control problem.
\end{rem}

The solution for $z(t)$ is found by solving the following minimization problem \cite{c20}:
\begin{equation} 
\begin{aligned} 
& \underset{z_0}{\text{minimize}}
& &  z_0^TP(0)z_0 , 
\label{minimization}
\end{aligned}
\end{equation}
where $P(0)$ is the value of $P(t)$ at $t_0$ and $P(t)$ is the solution to the following Riccati equation \cite{c20} with $P(t_f)=\mathbf{0}$:
    \begin{align}
        \dot{P}(t)&=-P(t)A_1(t)-A_1^{T}(t)P(t)-Q_1(t)\notag\\&+(P(t)B(t)+S(t))R^{-1}(t)(B^{T}(t)P(t)+S^{T}(t)).
      \label{eqn:secondary_riccati}
    \end{align}

\subsection{Hard-constrained method}
In the hard-constrained method \cite{c10a}, (\ref{first_cond}) is considered to be strictly satisfied while (\ref{second_cond}) is considered to be satisfied only approximately. In this method, the left-hand side of (\ref{second_cond}) is minimized assuming (\ref{first_cond}) holds strictly. In \cite{c10a}, constraints are imposed on the control trajectories in the primary problem. We will not consider constraints on the control trajectories in this work. The IOC problem for the method is:
\begin{equation}
\begin{aligned}
& \underset{c}{\text{minimize}}
& &   \int_{t_0}^{t_f} \Vert \nabla_u H(t,x(t),u(t)) \Vert^2 dt \\
& \text{subject to}
& & \dot{p}(t)=-\nabla_x H(t,x(t),u(t))\\
& & & p(t_f)=\mathbf{0}.
\label{hard_optimization}
\end{aligned}
\end{equation}
In \cite{c10a}, a variable, $L(t)$, is considered such that $p(t)=L(t)c$. Thus, (\ref{hard_optimization}) is reformulated as:

\begin{equation}
\begin{aligned}
& \underset{c}{\text{minimize}}
& c^TWc\\
& \text{subject to}
& & \dot{L}(t)=-\nabla_x\phi^T(t,x(t),u(t))\\
& & &  \hspace{0.65cm}-\nabla_x f^T(t,x(t),u(t))L(t), \\
& & & L(t_f)=\mathbf{0}.
\label{ioc_hard_optimization}
\end{aligned}
\end{equation}

Here, $W=\int_{t_0}^{t_f} W_1^T(t)W_1(t) dt$ and $W_1(t)=\nabla_u\phi^T(t,x(t),u(t))+\nabla_u f^T(t,x(t),u(t))L(t)$. $W_1(t)$ is obtained from $\nabla_u H(t,x(t),u(t))$. 
\begin{rem}
    For the hard-constrained method, (\ref{ioc_hard_optimization}) is the secondary optimization problem.
\end{rem}

\textbf{Note:} In \cite{c10a}, the inverse differential game problem is considered. In this work, the analysis is restricted to a single player, that is, the IOC problem.

\subsection{Motivation for the analysis of solvability of the minimum principle based IOC methods} 
It was stated in \cite{c20} that it is not clear when the soft-constrained method would fail. In \cite{c10a}, it was stated that addressing the selection of suitable dynamics was beyond the scope of their work. This work provides some answers to an open question about the kinds of optimal trajectories, the initial conditions of the closed-loop system and the system dynamics for the original optimal control problem, for which both the minimum principle based IOC methods will find the true solution. The paper provides cases where the soft-constrained method is able to solve the IOC problem for certain kinds of trajectories, which fails for the hard-constrained method.

\section{Solvability of the IOC problem for soft and hard-constrained methods for a general nonlinear primary problem}
In this section, the groundwork will be laid for the analysis of the solvability of the IOC problem using both the soft and the hard-constrained IOC methods for a general nonlinear primary problem.

\subsection{Solvability of the IOC problem for the soft-constrained method}
Here, we will look into the solvability conditions for the soft-constrained method.
\subsubsection{A more convenient arrangement of the Riccati equation}
We now rearrange the terms of the Riccati equation (\ref{eqn:secondary_riccati}) for carrying out the analysis of the solvability of the IOC problem for the soft-constrained method.
\begin{align}
    \dot{P}(t)&=-A^T(t)P(t)-P(t)A(t)-Q(t)\notag\\&+P(t)B(t)R^{-1}(t)B^T(t)P(t),
    \label{final_riccati}
\end{align}
where
\begin{align}
A(t)&=A_1(t)-B(t)R^{-1}(t)S^T(t)=-\begin{bmatrix}
             \mathbf{0}\\
             I
\end{bmatrix}S^T(t)\notag\\&=-\begin{bmatrix}
             \mathbf{0}\\
             I
\end{bmatrix}G^T(t)F(t)= -\begin{bmatrix}
\mathbf{0} & \mathbf{0}\\
\nabla_x\phi\vert_{(x,u)}^T & \nabla_xf\vert_{(x,u)}^T \\
\end{bmatrix}\label{eqn:A}
\end{align}
and
\begin{align}
Q(t)&=Q_1(t)-S(t)R^{-1}(t)S^T(t)=F^T(t)\begin{bmatrix}
             \mathbf{0} & \mathbf{0}\\
             \mathbf{0} & I
\end{bmatrix}F(t)\notag\\&=\begin{bmatrix}
\nabla_u\phi\vert_{(x,u)} \\ \nabla_uf\vert_{(x,u)} \\
\end{bmatrix}\begin{bmatrix}
\nabla_u\phi\vert_{(x,u)}^T & \nabla_uf\vert_{(x,u)}^T \\
\end{bmatrix}=C^T(t)C(t),
\label{eqn:Q}
\end{align}
for
\begin{equation}\label{equ:Ct}
C(t)=\begin{bmatrix}
\nabla_u\phi\vert_{(x,u)}^T & \nabla_uf\vert_{(x,u)}^T \\
\end{bmatrix}.
\end{equation}
$A(t)$, $B(t)$ and $C(t)$ from equations (\ref{eqn:A}) and (\ref{eqn:Q}) will be used for analysing the solvability of the IOC problem.
\subsubsection{Time-varying observability of the secondary problem}
It is now logical to state a remark based on the Riccati solution of the secondary optimal control problem.
\begin{rem} \label{rem:solvability}
  The solution  of equation (\ref{minimization}) exists and is unique if $P(0)$ is unique and $P(0)$ is positive definite \cite[Theorem 2]{c10a}.  Given that generally some of the parameters in the vector $\bar{c} \in \mathbb{R}^k$ are known, this condition reduces to the principal minor of $P(0)$ which excludes the columns for the known parameters \cite[Theorem 2]{c10a}. An alternative compact formulation is: Given a vector $\bar{c} \in \mathbb{R}^k$ populated only by non-zero elements, where the elements in $c$ are known. It is sufficient $N_s^T P(0) N_s$ to be full rank for \begin{equation}N_s=\ker\left(\text{diag}\left(\begin{bmatrix} \bar{c} \\ \mathbf{0}_{n \times 1} \end{bmatrix}\right)\right). \label{equ:Ns}\end{equation}
\end{rem}

\vspace{0.2cm}
Following Remark \ref{rem:solvability}, to gain deeper insight into the solvability of the IOC problem, it is necessary to further analyse the solvability of the Riccati equation (\ref{final_riccati}) equivalent to equation  (\ref{eqn:secondary_riccati}). If the linear time-varying system $(C(t),A(t))$ is observable, then the Riccati equation has a global unique solution and the solution, $P(t)$, is invertible \cite{c22}. This leads to a sufficient result for solvability of the IOC problem, which is summarized below: 
\begin{thm}\label{thm:1}  It is assumed that $C(t)$, $A(t)$ and $B(t)$ exist and  are bounded within the interval $[t_0,\,t_f]$. There is an interval $[t_0,\,T]$ for $0<T\leq t_f$, where the pair $(C(t),\,A(t))$ is observable. Under this condition, $P(0)$ is positive definite and a unique solution to the IOC problem exists.
\end{thm}
{\emph{Proof:}}  As $C(t)$, $A(t)$, $B(t)$ exist and are bounded, \cite[Lemma 1 and 2]{c22} imply that a unique and positive semi-definite solution $P(t)$ to the time-varying Riccati equation (\ref{final_riccati}) exists, given $P(t_f)$ is positive semidefinite and bounded. \cite[Theorem 1]{c22} implies that observability of $(C(t),A(t))$ in a finite interval $[0,T]$ guarantees a positive definite $P(0)$. From Remark \ref{rem:solvability}, this implies unique solvability of the IOC problem. {\hfill $\blacksquare$}

The observability for the linear time-varying system of $(C(t), A(t))$ of
order $l=k+n$ is best analysed using the observability matrix as given in \cite[Equation 4.107, p.117]{c23}:
\begin{equation}
    Q_o(t)=\begin{bmatrix}
    C^T(t)& \Delta_oC^T(t) & \ldots & \Delta_o^{l-1}C^T(t)
    \end{bmatrix},
    \label{equ:Qo}
\end{equation}
where $\Delta_o\equiv A^T(t)+\frac{d}{dt}$. Here, $A^T(t)=\begin{bmatrix}
\mathbf{0} & -\nabla_x\phi\vert_{(x,u)}\\
\mathbf{0} & -\nabla_xf\vert_{(x,u)} 
\end{bmatrix}$. Given the definition of  $Q_o(t)$ above, a less conservative sufficient theorem for solvability  of the IOC problem is:
\begin{thm} \label{thm:2}  It is assumed that $C(t)$, $A(t)$ , $B(t)$ and $Q_o(t)$ exist and  are bounded within the interval $[t_0,\,t_f]$. The matrix $N_s$ in (\ref{equ:Ns}) is defined by the known parameters in the vector $c$. There is an interval $[t_0,\,T]$ for $0<T\leq t_f$, where $Q_p=N_s^T Q_o(t) Q_o(t)^T N_s $, is full rank almost everywhere. Under this condition $N_s^T P(0) N_s$ is positive definite and in $[t_0,\,T]$, a unique solution to the IOC exists.

\noindent {\emph{Proof:}}  As in Theorem \ref{thm:1}, a unique and positive semi-definite solution, $P(t)$, of the time-varying Riccati equation (\ref{final_riccati}) exists, given $P(t_f)$ is positive semidefinite and bounded. From \cite[Equation (3.5)]{c10b}, $N_s^T P(0) N_s\geq N_s^T \int_{t_0}^{T} \Phi^T(s,t_0) C^T(t) C(t) \Phi(s,t_0) ds N_s$, where $\Phi(s,t_0)$ is the state transition matrix of $A(t)$. The latter matrix is positive definite since $Q_p=N_s^T Q_o(t) Q_o(t)^T N_s $ is full rank almost everywhere in $[t_0,\,T]$ \cite[Section 4.2]{c23}.  From Remark \ref{rem:solvability}, this implies unique solvability of the IOC. 
{\hfill $\blacksquare$   }\end{thm}

The observability matrix in (\ref{equ:Qo}) for the pair, $(C(t),A(t))$, obtained from (\ref{eqn:A}) and (\ref{equ:Ct}), is thus, $Q_o=\begin{bmatrix}
Q_{o1} & Q_{o2} & \ldots & Q_{oi} & \ldots & Q_{ol}
\end{bmatrix}$, where $Q_{oi}=\Delta_o^{i-1}C^T(t)$,~ ($i=1,\ldots,l$) and $Q_{o1}=\begin{bmatrix}
    (\nabla_u\phi)^T & (\nabla_uf)^T
    \end{bmatrix}^T$, $Q_{o2}=\begin{bmatrix}
    (\begin{smallmatrix}-\nabla_x\phi\nabla_uf+\frac{d}{dt}(\nabla_u\phi)\end{smallmatrix})^T & (\begin{smallmatrix}-\nabla_xf\nabla_uf+\frac{d}{dt}(\nabla_uf)\end{smallmatrix})^T
    \end{bmatrix}^T$.
This can be continued in a similar way up to the $l^{\text{th}}=(k+n)^{\text{th}}$ element of $Q_o$ as $Q_{o(n+k)}=\left(\begin{bmatrix}
\mathbf{0} & -\nabla_x\phi\vert_{(x,u)}\\
\mathbf{0} & -\nabla_xf\vert_{(x,u)} 
\end{bmatrix}+\frac{d}{dt}\right)^{n+k-1}\begin{bmatrix}
   \nabla_u\phi\vert_{(x,u)} \\
   \nabla_uf\vert_{(x,u)} 
\end{bmatrix}$.
\begin{rem}
It can be seen that the last $n$ rows of the observability matrix, $Q_o$, is the controllability matrix of the original optimal control problem, while the remaining rows are formed using the trajectories and system dynamics.
\label{controllability}
\end{rem}

\subsection{Solvability of the IOC problem for the hard-constrained method}
The solvability of the hard-constrained method will be analysed here. 
\begin{rem} \label{rem:solvability_hard}
  The solution of (\ref{ioc_hard_optimization}) exists and is unique if $W$ is unique and positive definite \cite[Theorem 3]{c10a}. Given that generally some of the parameters in the vector $c$ are known, this condition reduces to the principal minor of $W$ which excludes the columns for the known parameters \cite[Theorem 3]{c10a}. An alternative compact formulation is: Given a vector $\bar{c} \in \mathbb{R}^k$ populated only by non-zero elements, where the elements in $c$ are known. It is sufficient $N_h^T W N_h$ to be full rank for \begin{equation}N_h=\ker\left(\text{diag}\left( \bar{c}  \right)\right). \label{equ:Nh}\end{equation}
\end{rem}
We now look at a theorem for the solvability of the IOC problem for the hard-constrained method.
\begin{thm}
The IOC problem in (\ref{ioc_hard_optimization}) will have a unique solution if $N_h^TWN_h$ is non-singular.

\emph{Proof:} The proof can be found in \cite[Theorem 3]{c10a}.
\hfill$\blacksquare$
\end{thm}
We would now consider different kinds of optimal closed-loop trajectories and initial conditions for the primary problem, to analyse the rank of $N_h^TWN_h$ and extract it from $W$. For this, it makes sense to analyse $W_1^T(t)W_1(t)$.
\begin{align}
    &W_1^T(t)W_1(t)=\nabla_u\phi(t)\nabla_u\phi^T(t)+\nabla_u\phi(t)\nabla_u f^T(t)L(t)\notag\\&+L^T(t)\nabla_u f(t)\nabla_u\phi^T(t)+L^T(t)\nabla_u f(t)\nabla_u f^T(t)L(t).
    \label{w1}
\end{align}
In (\ref{w1}), $L(t)$ is found by solving the differential equation in (\ref{ioc_hard_optimization}) backwards in time, with terminal boundary condition, $L(t_f)=\mathbf{0}$. Thus, we arrive at
\begin{align}
L(t)&=\int_{t}^{t_f}\psi(t,\tau)\nabla_x\phi^T(\tau)d\tau,
\label{l}
\end{align}
where $\psi(t,\tau)$ is the state transition matrix of $-\nabla_x f^T(t,x(t),u(t))$.  We then substitute this expression of $L(t)$ obtained in (\ref{l}) into (\ref{w1}) to arrive at (\ref{w1w1}) shown at the bottom of page \pageref{w1w1}. From (\ref{w1w1}), if $N_h^TWN_h$ is non-singular,  where $W=\int_{t_0}^{t_f} W_1^T(t)W_1(t) dt$, the IOC problem will have a unique solution \cite[Theorem 3]{c10a}.

\begin{figure*}[b]
\hrulefill
 \begin{equation}
 \begin{aligned} \label{w1w1} 
W_1(t)^TW_1(t)&=\nabla_u\phi(t)\nabla_u\phi^T(t)+\nabla_u\phi(t)\nabla_u f^T(t)\Biggl(\int_{t}^{t_f}\psi(t,\tau)\nabla_x\phi^T(\tau)d\tau\Biggr)+\Biggl(\int_{t}^{t_f}\nabla_x\phi(\tau)\psi^T(t,\tau)d\tau\Biggr)\nabla_u f(t)\times\\&\nabla_u\phi^T(t)+\Biggl(\int_{t}^{t_f}\nabla_x\phi(\tau)\psi^T(t,\tau)d\tau\Biggr)\nabla_u f(t)\nabla_u f^T(t)\Biggl(\int_{t}^{t_f}\psi(t,\tau)\nabla_x\phi^T(\tau)d\tau\Biggr)   
\end{aligned}
\end{equation}
\end{figure*}

\section{Solvability of the IOC problem for soft and hard-constrained methods for an infinite horizon Linear Time Invariant primary problem}

This sections deals with the solvability of both the IOC methods when the primary problem is an infinite horizon LTI system with a single scalar input $u\in \mathbb{R}$. 

\begin{equation}
\begin{aligned}
& \underset{u(t)}{\text{minimize}}
& & \int_{t_0}^{\infty} (x^T(t)Dx(t)+u^T(t)Eu(t)) dt \\
& \text{subject to}
& & \dot{x}(t)=Mx(t)+Nu(t)\\
& & & x(t_0)=x_0.
\end{aligned}
\label{primary_optimal_problem}
\end{equation}

\begin{assu}
It is assumed that the matrices $D$ and $E$ are diagonal and $E$ is the only known coefficient.
\end{assu}
Note that this assumption implies that the elements of the weight, $D$, are non-negative, and $E$ is positive. For this infinite horizon optimal control problem, the control input is $u(t)=\Theta x(t)$, where $\Theta=-E^{-1}N^T\Pi$ is the control gain for the infinite-horizon LQR problem. $\Pi$ is found by solving the following Algebraic Riccati Equation, $M^T\Pi+\Pi M+D-\Pi NE^{-1}N^T\Pi=0$. $u(t)$ can also be written as $u(t)=-E^{-1}N^Tp(t)$, where $p(t)=\Pi x(t)$ and $p(\infty)=0$. The cost function can be written as:
\begin{equation}
\begin{aligned}
   & \int_{t_0}^{\infty} c^T\phi(t,x(t),u(t)) dt, 
\end{aligned}
\end{equation}
where $c=\begin{bmatrix} c_1 & c_2 & \ldots c_{n+1}\end{bmatrix}^T$ and $\phi=\begin{bmatrix} x_1^2(t) & x_2^2(t) & \ldots & x_n^2(t) &  u^2(t) \end{bmatrix}$.
\begin{rem}
The goal of the IOC problem is to recover the weights, $c$, of the cost function by measuring $x(t)$ and $u(t)$.
\end{rem}

\subsection{Solvability of the IOC problem using the soft-constrained method for an infinite horizon LTI Primary Optimal Control Problem}
We will now consider the solvability of the IOC problem using the soft-constrained method for an infinite horizon LTI primary problem. In this case, for (\ref{eqn:secondary_lqr}), $z=\begin{bmatrix}c_1 & c_2 & \ldots & c_{n+1} & p_1 & p_2 & \ldots p_n\end{bmatrix}^T$, $v=\begin{bmatrix}\dot{p}_1 &\dot{p}_2 & \ldots &\dot{p}_n\end{bmatrix}^T$, $\nabla_xf=M$, $\nabla_uf=N$, $\nabla_u\phi=\begin{bmatrix} 0 & 0 & \dots & 2u \end{bmatrix}^T$ and $\nabla_x\phi=\hbox{diag}(2x_1,2x_2, \ldots, 2x_n,0)$. Hence, for the second-order linear system in (\ref{primary_optimal_problem}), 
$F=\begin{bmatrix}
 \nabla_x\phi\vert_{(x,u)}^T & M^T \\
\nabla_u\phi\vert_{(x,u)}^T & N^T \end{bmatrix}$ and $G=\begin{bmatrix} I\\\mathbf{0} \end{bmatrix}$. From the observability matrix, $Q_o$, of the secondary problem for a primary infinite horizon LTI problem, the $i^{\text{th}}$ element of $Q_o$ for the general LTI system in (\ref{primary_optimal_problem}) is:
\begin{equation}
 Q_{o\,i}=  - (-1)^i \begin{bmatrix}
    (\nabla_x\phi) M^{i-2}N\\ M^{i-1}N
    \end{bmatrix}+ \frac{d}{dt} Q_{o\,(i-1)},~i=1, \ldots l.
\end{equation}
For instance, for a second-order linear system in (\ref{primary_optimal_problem}), (\ref{equ:Qo}) can be written as (\ref{equ:Qo_second_gen}) as shown at the bottom of page \pageref{equ:Qo_second_gen}.
\begin{figure*}[b]
\hrulefill
\begin{equation}
\begin{aligned}
\label{equ:Qo_second_gen}
    Q_o=\begin{bmatrix}
  (\nabla_u\phi)^T & N^T \\
  (-(\nabla_x\phi)N+\frac{d}{dt}(\nabla_u\phi))^T & -(MN)^T \\
  ((\nabla_x\phi)MN-(\frac{d}{dt}(\nabla_x\phi))N+\frac{d^2}{dt^2}(\nabla_u\phi))^T & (M^2N)^T \\ (-(\nabla_x\phi)M^2N+\{\frac{d}{dt}(\nabla_x\phi) \}MN-\{\frac{d^2}{dt^2}(\nabla_x\phi)\}N+\frac{d^3}{dt^3}(\nabla_u\phi))^T & (-M^3N)^T \\ ((\nabla_x\phi)M^3N-\{\frac{d}{dt}(\nabla_x\phi)\}M^2N +\{\frac{d^2}{dt^2}(\nabla_x\phi)\}MN-\{\frac{d^3}{dt^3}(\nabla_x\phi)\}N+\frac{d^4}{dt^4}(\nabla_u\phi))^T & (M^4N)^T
    \end{bmatrix}^T
    \end{aligned}
\end{equation}
\hrulefill
\end{figure*}

Note again that from \textit{Remark \ref{controllability}}, it can be understood that the lower $n$ rows of the observability matrix, $Q_o$, is the controllability matrix of the primary LTI problem, while the upper $k$ rows of $Q_o$ are dependent on the trajectories, $x(t)$ and $u(t)$, and the system matrices, $M$ and $N$.

\subsubsection{Analysis for different types of trajectories using soft-constrained method for an infinite horizon second-order LTI primary problem}
It is known that second-order systems can be seen as the building blocks of higher-order systems. Hence, it is logical to consider the solvability of the IOC problem for second-order systems.  We will consider different kinds of trajectories and initial conditions from an infinite horizon linear second-order primary optimal control problem.  In this case, the time-varying observability matrix in (\ref{equ:Qo}) is $Q_o\in \mathbb{R}^{5\times 5}$ and $z(t)=\begin{bmatrix}
  c_1 & c_2 &c_3 & p_1(t) & p_2(t)  \end{bmatrix}^T$. Here, $c_3$ is assumed to be known.\\
\textbf{Case 1}: We will now consider trajectories obtained from a second-order infinite-horizon over-damped closed-loop primary problem initialised at a single real mode.
\begin{cor} 
    For an over damped closed-loop primary LTI problem (\ref{primary_optimal_problem}) initialised at a single real mode, using the soft-constrained method, the IOC problem (\ref{minimization}) will not converge to the true solution.

{\emph{Proof:}}  Let us assume $(M+N\Theta)$ has the real eigenvalues $\lambda_i$ with right eigenvector $V_i$ and we consider one of the eigenvalues and limit the dynamic closed-loop system to that mode, i.e. $x(t_0)=V_1$ and $x(t)=V_1 e^{\lambda_1 t}$, where $\lambda_1 \in \eigLambda(M+N\Theta)$. This has the following implications on (\ref{equ:Qo}) as shown in (\ref{Qo_over}) at the bottom of page \pageref{Qo_over}.
\begin{figure*}[b]
\hrulefill
\begin{equation}
\begin{aligned} \label{Qo_over}
    Q_o=\begin{bmatrix}
    (\nabla_u\phi)^T & N^T \\
    (\nabla_x\phi(-N)+\lambda_1\nabla_u\phi)^T & (-MN)^T \\
    (\nabla_x\phi(MN-\lambda_1N)+\lambda_1^2\nabla_u\phi)^T & (M^2N)^T\\
    (\nabla_x\phi(-M^2N+\lambda_1MN-\lambda_1^2M)+\lambda_1^3\phi_u)^T & (-M^3N)^T\\
    (\nabla_x\phi(M^3N-\lambda_1M^2N+\lambda_1^2MN-\lambda_1^3N)+\lambda_1^4\nabla_u\phi)^T & (M^4N)^T
    \end{bmatrix}^T
    \end{aligned}
\end{equation}
\end{figure*}
Writing $\nabla_x\phi=2\begin{bmatrix}
\text{diag}(V_1e^{\lambda_1t})\\\mathbf{0}_{1\times 2}
\end{bmatrix}$ and $\nabla_u\phi=2\begin{bmatrix} 0 & 0 & \Theta V_1e^{\lambda_1t}\end{bmatrix}^T$, (\ref{Qo_over}) can then be re-written as $Q_o=(Q_{om}  Q_{oo})^T$, where $Q_{om}=\begin{bmatrix}
2\text{diag}(V_1e^{\lambda_1t}) & \mathbf{0}_{2\times 1} & \mathbf{0}_{2\times 2}\\\mathbf{0} _{1\times 2} & 2\Theta V_1e^{\lambda_1t} & \mathbf{0}_{1\times 2}\\\mathbf{0}_{2\times 2} & \mathbf{0}_{2\times 1} & I_{2\times 2}
\end{bmatrix}$ and $Q_{oo}=\begin{bmatrix}
\mathbf{0}_{2\times 1} & -N & MN-\lambda_1N & \begin{smallmatrix}(-M^2N\\+\lambda_1MN\\-\lambda_1^2N)\end{smallmatrix} & \begin{smallmatrix}(M^3N\\-\lambda_1M^2N\\+\lambda_1^2MN-\lambda_1^3N)\end{smallmatrix}\\1 & \lambda_1 & \lambda_1^2 & \lambda_1^3 & \lambda_1^4\\N & -MN & M^2N & -M^3N & M^4N
\end{bmatrix}$. If $V_1=\begin{bmatrix}
    v_{11} & v_{12}
\end{bmatrix}^T$ and  if $v_{11}\neq 0$, $v_{12}\neq 0$ and $\Theta V_1 \neq 0$, $Q_{om}$ is full rank. Hence, the rank of $Q_o$ depends on the rank of $Q_{oo}$. Let, $\beta=\begin{bmatrix}
    \mathbf{0}_{2\times 1}\\ 1\\N
    \end{bmatrix}$ and $\alpha=\begin{bmatrix}
    \vert \lambda_1\vert I_{2\times 2} & \mathbf{0}_{2\times 1} & I_{2\times 2}\\\mathbf{0}_{1\times 2} & \vert\lambda_1\vert & \mathbf{0}_{1\times 2}\\\mathbf{0}_{2\times 2} & \mathbf{0}_{2\times 1} & M
    \end{bmatrix}$. $Q_{oo}$ can be seen as the controllability matrix,
    $Q_{oo}=\begin{bmatrix}
    \beta & -\alpha\beta & \alpha^2\beta & -\alpha^3\beta & \alpha^4\beta
    \end{bmatrix}$, for the pair $(\alpha, \beta)$. We used the PBH test to determine the rank with $s=\vert\lambda_1\vert$: $\begin{bmatrix}
        sI-\alpha & \beta
        \end{bmatrix}=\begin{bmatrix}
        \mathbf{0}_{2\times 2} & \mathbf{0}_{2\times 1} & -I_{2\times 2} & \mathbf{0}_{2\times 1}\\\mathbf{0}_{1\times 2} & 0 & 0 & 1\\\mathbf{0}_{2\times 2} & \mathbf{0}_{2\times 1} & \vert\lambda_1\vert I_{2\times 2}-M & N
        \end{bmatrix}$. This matrix will have a maximum rank of 3 due to two zero columns in a $5\times 5$ matrix. Hence, the IOC problem will not converge to the true solution.
        {\hfill $\blacksquare$}
    \end{cor}
The analysis for a critically damped system initialised at a single real mode is the same as this case, and will be omitted. 

\textbf{Case 2:} We will now consider the case when the trajectories are obtained from an infinite horizon second-order under-damped closed-loop primary problem. Let us assume that the closed-loop system matrix of the primary problem is $\bar{M}=\begin{bmatrix}
    \sigma & \omega \\-\omega & \sigma
    \end{bmatrix}$, the initial states are $\begin{bmatrix}
    x_{1,0} & x_{2,0} \end{bmatrix}^T$ and the complex eigenvalues of the closed-loop system are $\sigma \pm j\omega$, where $\sigma$ is real and negative. Thus, $Q_o$ in (\ref{equ:Qo}) now becomes
    \begin{align}
    Q_o&=\begin{bmatrix}
    2e^{\sigma t}(Q_{cc}\cos(\omega t)+Q_{ss}\sin(\omega t)) & q_c
    \end{bmatrix}^T,
    \label{Qo_under}
   \end{align} 
where
\begin{equation}
q_c=\begin{bmatrix}
N^T & (-MN)^T & (M^2N)^T & (-M^3N)^T & (M^4N)^T
\end{bmatrix}^T
\end{equation}
and, $Q_{cc}=\begin{bmatrix}
    Q_{oc1} & Q_{oc2} & Q_{oc3} & Q_{oc4} & Q_{oc5}
    \end{bmatrix}^T$, where $Q_{oci}=\begin{bmatrix}
    Q_{ocai} & Q_{ocki} \end{bmatrix}$ for $i=1,..,5$ and $Q_{ss}=\begin{bmatrix}
    Q_{os1} & Q_{os2} & Q_{os3} & Q_{os4} & Q_{os5}
    \end{bmatrix}^T$, where $Q_{osi}=\begin{bmatrix}
    Q_{osai} & Q_{oski}
    \end{bmatrix}$ for $i=1,..,5$. $Q_{oca1}=Q_{osa1}=\begin{bmatrix}
      0 & 0
    \end{bmatrix}^T$ and
    \begin{equation}
    Q_{ocai}=    \sum_{l=2}^{i} (-1)^{i-l+1}\hbox{diag}\left(\bar{M}^{l-2}\begin{bmatrix}
    x_{1,0}\\x_{2,0}
    \end{bmatrix}\right)M^{i-l}N
    \end{equation}
    \begin{align}
    Q_{osai}&= \sum_{l=2}^{i} (-1)^{i-l+1}\hbox{diag}\left(\begin{bmatrix}
    1 & 0\\0 & -1
    \end{bmatrix}(\bar{M}^T)^{l-2}\begin{bmatrix}
    x_{2,0}\\x_{1,0}
    \end{bmatrix}\right)\times\notag\\&M^{i-l}N,
    \end{align}
    for $i=2, \ldots, 5$. Moreover,
    \begin{equation}
    Q_{ocki}=\begin{bmatrix}
    \theta_1  & \theta_2
    \end{bmatrix}\bar{M}^{i-1}\begin{bmatrix}
     x_{1,0}\\ x_{2,0}
    \end{bmatrix}
    \end{equation}
    \vspace{-0.2cm}
\begin{equation}
Q_{oski}=\begin{bmatrix}
    \theta_1 & \theta_2
    \end{bmatrix}\begin{bmatrix}
    1 & 0\\0 & -1
    \end{bmatrix}(\bar{M}^T)^{i-1}\begin{bmatrix}
    x_{2,0} \\ x_{1,0}
    \end{bmatrix},
    \end{equation}
for $i=1, \ldots, 5$. To find the rank of the observability matrix, $Q_o$, for this case, we carried out a numerical simulation (in Section V) and found that for a primary underdamped problem (\ref{primary_optimal_problem}) using the soft-constrained method, the IOC problem (\ref{minimization}) converged to the true solution.
    \begin{rem}
        It can be understood that $N_s^TQ_o$ needs to be full rank  in an interval $[t_0,T]$ for the IOC problem to be solvable (Theorem \ref{thm:2}). In this case, it can be of interest to observe the two terms in (\ref{Qo_under}) associated with the sine and cosine terms.
    \end{rem}

\subsection{Solvability of the IOC problem using the hard-constrained method for an Infinite Horizon LTI Primary Problem}
We will now consider the primary problem given in (\ref{primary_optimal_problem}). Thus, for this problem, $L(t)$ in (\ref{l}) becomes
\begin{align}
    L(t)&=\int_{t}^{\infty}e^{-M^T(t-\tau)}\nabla_x\phi^T(\tau)d\tau.
    \label{l_linear}
\end{align}
Substituting $L(t)$ from (\ref{l_linear}), $W$ in (\ref{ioc_hard_optimization}) can be written as (\ref{w_general}) shown at the bottom of page \pageref{w_general}.

\begin{figure*}[b]
\hrulefill
\begin{equation}
\begin{aligned}
    W&=4\int_{t_0}^{\infty}\begin{bmatrix}
        \begin{smallmatrix}\{\int_{t}^{\infty}\text{diag}(x(\tau))e^{-M(t-\tau)}d\tau\}NN^T\{\int_{t}^{\infty}e^{-M^T(t-\tau)}\text{diag}(x(\tau))d\tau\} \end{smallmatrix} & \begin{smallmatrix}\{\int_{t}^{\infty} \text{diag}(x(\tau))e^{-M(t-\tau)}d\tau\}N\Theta x(t)\end{smallmatrix} \\ \begin{smallmatrix}\Theta x(t)N^T\{\int_{t}^{\infty}e^{-M^T(t-\tau)}\text{diag}(x(\tau))d\tau\}\end{smallmatrix} & (\Theta x(t))^2
    \end{bmatrix} dt
    \label{w_general}
\end{aligned}
\end{equation}
\hrulefill
\end{figure*}

\subsubsection{Analysis for different kinds of trajectories of second-order primary problems using hard-constrained methods}
In this section, the solvability of the hard-constrained IOC method for different kinds of optimal trajectories and initial conditions from an infinite horizon second-order LTI closed-loop primary problem will be analysed. The focus is on second-order systems. Here, $W\in\mathbb{R}^{3\times 3}$. The element $c_3$ of the vector $c=\begin{bmatrix}
    c_1 & c_2 & c_3
\end{bmatrix}^T$ is assumed to be known.

\textbf{Case 1}: We will consider trajectories obtained from an over-damped infinite-horizon closed-loop primary problem initialised at a single real mode. We will use the same problem here which was used for the over-damped case in the soft-constrained method. $W$ in (\ref{w_general}) can be written as $W=4\int_{t_0}^{\infty}W_a(t) dt$, where $W_a(t)$ is obtained from (\ref{w_general}) by substituting $\hbox{diag}(x(\tau))=\hbox{diag}(V_1)e^{\lambda_1 \tau}$. Next, we can write $W_a(t)=W_{22}W_{2a}(t)W_{22}$, where $W_{22}=\begin{bmatrix}
        \text{diag}(V_1) & \mathbf{0}_{2\times 1} \\ \mathbf{0}_{1\times 2} & \Theta V_1
    \end{bmatrix}$ and $W_{2a}(t)=W_{2a1}(t)W_{2a1}^T(t)$, where $W_{2a1}(t)$ is:
\begin{equation}\label{w2a1}
W_{2a1}(t)=\begin{bmatrix}
       e^{-Mt}\{\int_{t}^{\infty} e^{(\lambda_1 I+M)\tau} d\tau\}N\\e^{\lambda_1 t}
    \end{bmatrix}.
\end{equation}
If $V_1=\begin{bmatrix}
    v_{11} & v_{12}
\end{bmatrix}^T$ and $v_{11}\neq 0$, $v_{12}\neq 0$ and $\Theta V_1 \neq 0$, $W_{22}$ is full rank. The rank of $W$, thus, boils down to the rank of $\int_{t_0}^{\infty}W_{2a}(t)dt$. Let, $\bar{W}_{a1}(t)=e^{-Mt}\{\int_{t}^{\infty} e^{(\lambda_1 I+M)\tau} d\tau\}N$. The IOC problem has a unique solution if $W_{\psi}=\int_{t_0}^{\infty}\bar{W}_a(t)dt=\int_{t_0}^{\infty}\bar{W}_{a1}(t)\bar{W}_{a1}^T(t)dt$ is full rank. The following Lemma is needed for analysis of $W_{\psi}$.

\begin{lem}
    Consider a matrix, $Y=\int_{t_0}^{\infty}y(t)y^T(t)dt$, where, $y(t)\in \mathbb{R}^n$ is continuously differentiable. $Y$ is full rank, if and only if there exists a set of points, $t_0,t_1,\ldots, t_i,\ldots,t_n\leq \lim_{T_1\to\infty}T_1$ for all non-zero vectors $q$, such that
\begin{equation}
    q^T\underbrace{\begin{bmatrix}
        y(t_0) & y(t_1) & \ldots y(t_i) & \ldots & y(t_n)
    \end{bmatrix}}_{Y_a}\neq 0.
    \label{y_eqn}
\end{equation}
\label{rank_lemma}
\end{lem}
 The proof can be easily shown following the arguments of \cite[p.610]{c25}. It is possible to use Lemma \ref{rank_lemma} for finding the rank of $W_{\psi}$ by defining
\begin{align}
    &y(t_i)=e^{-Mt_i}\{\int_{t_i}^{\infty} e^{(\lambda_1 I+M)\tau} d\tau\}N\notag\\&=\lim_{T_1 \to \infty}e^{-Mt_i}\{e^{(\lambda_1I+M) T_1}-e^{(\lambda_1I+M)t_i}\}(\lambda_1I+M)^{-1}N.
    \label{y_ti}
\end{align}
We need to find the rank of $Y_a$ now using (\ref{y_ti}). We will consider two scenarios for this analysis: (i) $\text{Re}(\bar{\lambda})<0$, is satisfied for any $\bar{\lambda} \in \eigLambda(\lambda_1 I+M)$ (ii) $\text{Re}(\bar{\lambda})>0$, is satisfied for some $\bar{\lambda} \in \eigLambda(\lambda_1 I+M)$. Here, $\lambda_1\in \eigLambda(M+N\Theta)$. \\
\textbf{\textit{Scenario 1}}: We will consider the case when $\text{Re}(\bar{\lambda})<0$, is satisfied for any $\bar{\lambda} \in \eigLambda(\lambda_1 I+M)$, where $\lambda_1\in \eigLambda(M+N\Theta)$. 
\begin{lem}
    For an over-damped closed-loop primary LTI problem (\ref{primary_optimal_problem}) initialised at a single real mode, using the hard-constrained method, the IOC problem (\ref{ioc_hard_optimization}) will not converge to the true solution, when $\text{Re}(\bar{\lambda})<0$, for any $\bar{\lambda} \in \eigLambda(\lambda_1 I+M)$, where $\lambda_1\in \eigLambda(M+N\Theta)$.

{\emph{Proof:}} Here,
\begin{align}
y(t_i)=-e^{\lambda_1 t_i}(\lambda_1 I+M)^{-1}N.
\label{yti_1}
\end{align}
Thus, each column of $Y_a$ in (\ref{y_eqn}) is a linear combination of the column, $(\lambda_1 I+M)^{-1}N$. Thus, the rank of $W_{\psi}$ is 1 and the IOC problem will not converge to the true solution.
{\hfill $\blacksquare$}
\end{lem}
\textbf{\textit{Scenario 2}}: We will consider the case when $\text{Re}(\bar{\lambda})>0$, for some $\bar{\lambda} \in \eigLambda(\lambda_1 I+M)$, where $\lambda_1\in \eigLambda(M+N\Theta)$. 
\begin{lem}
    For an over-damped closed-loop primary LTI problem (\ref{primary_optimal_problem}) initialised at a single real mode, using the hard-constrained method, the IOC problem (\ref{ioc_hard_optimization}) will not converge, when $\text{Re}(\bar{\lambda})>0$, for some $\bar{\lambda} \in \eigLambda(\lambda_1 I+M)$, where $\lambda_1\in \eigLambda(M+N\Theta)$.

{\emph{Proof:}} Here, $\lim_{T_1\to\infty}e^{(\lambda_1I+M)T_1}=\infty$ and thus,  $\lim_{T_1\to\infty}\Vert y(t)y^T(t)\Vert=\infty$. Hence, $W_{\psi}$ and $W$ does not exist. The IOC problem (\ref{ioc_hard_optimization}) will not converge.
{\hfill $\blacksquare$}
\end{lem}

The analysis for a critically damped system initialised at a single real mode is the same as this case, and will be omitted. 

\textbf{Case 2}: We will consider trajectories obtained from an under-damped closed-loop primary problem. The closed-loop primary system matrix is $\bar{M}=\begin{bmatrix}
    \sigma & \omega \\-\omega & \sigma
\end{bmatrix}$.  In (\ref{w_general}), if the $(2\times 2)$ upper left sub-matrix, $W_b=\int_{t_0}^{\infty}\{\int_{t}^{\infty}\text{diag}(x(\tau))e^{-M(t-\tau)}d\tau\}NN^T\{\int_{t}^{\infty}e^{-M^T(t-\tau)}\\\text{diag}(x(\tau))d\tau\}\}dt$, is full rank, we will be able to recover the true weights of the cost function of the primary optimal control problem. We need to analyse $W_{b2}(t)=\{\int_{t}^{\infty}\text{diag}(x(\tau))e^{-M(t-\tau)}d\tau\}N$. For the under-damped case, $x(\tau)=e^{\sigma\tau}\begin{bmatrix}
    \begin{smallmatrix}x_{1,0}\text{cos}(\omega\tau)+x_{2,0}\text{sin}(\omega\tau)\end{smallmatrix} & \begin{smallmatrix}-x_{1,0}\text{sin}(\omega\tau)+x_{2,0}\text{cos}(\omega\tau)\end{smallmatrix}
\end{bmatrix}^T$, where $x_0=\begin{bmatrix}
    x_{1,0} & x_{2,0}
\end{bmatrix}^T$ is the initial value of $x$. Thus, using Lemma \ref{rank_lemma} for this case, $y(t_i)$ in (\ref{y_eqn}) is
\begin{align}
    &y(t_i)=\int_{t_i}^{\infty} e^{\sigma\tau}\begin{bmatrix}
      \begin{smallmatrix}\{x_{1,0}\text{cos}(\omega\tau)\\+x_{2,0}\text{sin}(\omega\tau)\}\end{smallmatrix}& 0\\
      0 & \begin{smallmatrix}\{-x_{1,0}\text{sin}(\omega\tau)\\+x_{2,0}\text{cos}(\omega\tau)\}\end{smallmatrix}
    \end{bmatrix}e^{-M(t_i-\tau)}d\tau N.
    \label{y_ti_under}
\end{align}
Writing $\text{cos}(\omega\tau)$ and $\text{sin}(\omega\tau)$ in Euler form, (\ref{y_ti_under}) becomes
\begin{equation}
    y(t_i)=\frac{1}{2}(W_{b3}(t_i)+W_{b4}(t_i)),
    \label{wb2}
\end{equation}
where 
\begin{align}
    W_{b3}(t_i)&=\lim_{T_1\to\infty}x_a e^{-Mt_i}(M+(\sigma+j\omega)I)^{-1}(e^{\{M+(\sigma+j\omega)I\}T_1}\notag\\&-e^{\{M+(\sigma+j\omega)I\}t_i})N\label{general_wb3}\\
    W_{b4}(t_i)&=\lim_{T_1\to\infty}x_a^* e^{-Mt_i}(M+(\sigma-j\omega)I)^{-1}(e^{\{M+(\sigma-j\omega)I\}T_1}\notag\\&-e^{\{M+(\sigma-j\omega)I\}t_i})N.
    \label{general_wb4}
\end{align}
Here, $x_a=\begin{bmatrix}
        x_{1,0}-x_{2,0}j & 0\\0& x_{2,0}+x_{1,0}j
    \end{bmatrix}$. (\ref{general_wb3}) and (\ref{general_wb4}) can be written as:
\begin{align}
    W_{b3}(t_i)&=\lim_{T_1\to\infty}x_a\{e^{-Mt_i}e^{(M+(\sigma+j\omega)I) T_1}(M\notag\\&+(\sigma+j\omega)I)^{-1}N-e^{(\sigma+j\omega) t_i}(M+(\sigma+j\omega)I)^{-1}N\}\label{modified_wb3}\\
    W_{b4}(t_i)&=\lim_{T_1\to\infty}x_a^*\{e^{-Mt_i}e^{(M+(\sigma-j\omega)I) T_1}(M\notag\\&+(\sigma-j\omega)I)^{-1}N-e^{(\sigma-j\omega) t_i}(M+(\sigma-j\omega)I)^{-1}N\}.
    \label{modified_wb4}
\end{align}

\textbf{\textit{Scenario 1}} ($\text{Re}\{\hat{\lambda}+\sigma \}<0$ for any $\hat{\lambda} \in \eigLambda(M)$): It can be seen that the following holds:
\begin{align}
    W_{b3}(t_i)&=-e^{\sigma  t_i}x_a(M+(\sigma+j\omega)I)^{-1}e^{j\omega  t_i}N\label{wb3_stable}\\
    W_{b4}(t_i)&=-e^{\sigma t_i}x_a^*(M+(\sigma-j\omega)I)^{-1}e^{-j\omega  t_i}N.
    \label{wb4_stable}
\end{align}
Let, $H=x_a(M+(\sigma+j\omega)I)^{-1}$. Let $H_r$ and $H_c$ be two matrices consisting  of the real and imaginary parts of the matrix, $H$. Thus, (\ref{wb2}) becomes:
\begin{align}
    y(t_i)&=e^{\sigma t_i}(H_c \hbox{sin}(\omega  t_i)-H_r \hbox{cos}(\omega t_i)) N.
    \label{wb2_stable}
\end{align}
Let, 
\begin{equation}
((1/\omega)(M+\sigma I)^2+\omega I)^{-1}N=\begin{bmatrix}
   \delta_1 & \delta_2 
\end{bmatrix}^T
\label{delta1delta2}
\end{equation}
\begin{equation}
(1/\omega)(M+\sigma I)\begin{bmatrix}
   \delta_1 & \delta_2 
\end{bmatrix}^T=\begin{bmatrix}
    \mu_1 & \mu_2
\end{bmatrix}^T.
\label{mu1mu2}
\end{equation}
It can be seen that
\begin{equation}
  H_rN=\text{diag}(x_{1,0},x_{2,0})\begin{bmatrix}
    \mu_1 & \mu_2
\end{bmatrix}^T
\label{hrn}
\end{equation}
\begin{equation}
  H_cN=\text{diag}(-x_{2,0},x_{1,0})\begin{bmatrix}
   \delta_1 & \delta_2
\end{bmatrix}^T.
\label{hcn}
\end{equation}
\begin{lem}
    For an under-damped closed-loop primary LTI problem (\ref{primary_optimal_problem}), if $\text{Re}\{\hat{\lambda}+\sigma \}<0$ for any eigenvalue, $\hat{\lambda} \in \eigLambda(M)$ and for $\delta_1$, $\delta_2$, $\mu_1$ and $\mu_2$ in (\ref{delta1delta2}) and (\ref{mu1mu2}), the hard-constrained method has the following characteristics: (i) If $\delta_1\mu_2\delta_2\mu_1>0$, the IOC problem (\ref{ioc_hard_optimization}) converges to the true solution, for any non-zero initial states. (ii)  If $\delta_1\mu_2\delta_2\mu_1\leq 0$, there exist initial states for which the the IOC problem (\ref{ioc_hard_optimization}) does not converge to the true solution.\\
    \emph{Proof:} This is easily proved by checking if $H_rN$ and $H_cN$ from (\ref{hrn}) and (\ref{hcn}) respectively are linearly independent.
{\hfill $\blacksquare$}
\end{lem}

\textbf{\textit{Scenario 2}} ($\text{Re}\{\hat{\lambda}+\sigma \}>0$ for some $\hat{\lambda} \in \eigLambda(M)$): If we test via root-locus analysis using Chang Letov equation \cite{c24, c25}, we can find that this case is possible only if $M$ has unstable real eigenvalues. 
\begin{lem}
    For an under-damped closed-loop primary LTI problem (\ref{primary_optimal_problem}), if the open-loop system, $M$, has unstable real eigenvalues of the primary problem such that $\text{Re}\{\hat{\lambda}+\sigma \}>0$, for some $\hat{\lambda} \in \eigLambda(M)$, then using the hard constrained method, the IOC problem (\ref{ioc_hard_optimization}) does not converge.

\emph{Proof:} When $\text{Re}\{\hat{\lambda}+\sigma \}>0$, for some $\hat{\lambda} \in \eigLambda(M)$, $\lim_{T_1\to\infty}e^{\{M+(\sigma+j\omega)I\}T_1}\rightarrow \infty$ in (\ref{general_wb3}) and (\ref{general_wb4}). Hence,  $\lim_{T_1\to\infty}\Vert y(t)y^T(t)\Vert=\infty$. 
Thus, the IOC problem (\ref{ioc_hard_optimization}) does not coverge for this scenario.
\hfill$\blacksquare$
\end{lem}

\section{Examples}
In this section, we will consider three examples.

\textbf{Example 1:} Let us consider the following example
\begin{equation}
    \dot{x}=\begin{bmatrix}
    0 &  -1\\6 & 5
    \end{bmatrix}x+\begin{bmatrix}
    0\\1
    \end{bmatrix}u.
    \label{example2_underdamped}
\end{equation}
The eigenvalues of the open-loop system matrix, $M$ are 2 and 3. The system is initialised along $x_0=\begin{bmatrix}
    1 & -3
\end{bmatrix}^T$. The weights of the cost function are chosen as $D=\begin{bmatrix}
    32 & 0\\0 & 2
\end{bmatrix}$ and $E=1$, which results in an under-damped closed-loop system. One of the eigenvalues of the open-loop system results in $\text{Re}(\hat{\lambda}+\sigma)>0$. (i) Using the soft-constrained method, the weights were recovered to be close to the true weights. The condition numbers of $Q_p$ was found to be small (implying $Q_p$ is full rank) indicating that the IOC problem is solvable by the soft-constrained method for under-damped trajectories. (ii) Using the hard-constrained IOC method, the recovered weights were not close to the true weights. The elements of the matrix, $Y_a$, could not be computed. Hence, the optimization problem in (\ref{ioc_hard_optimization}) does not converge for this example. 

\textbf{Example 2}: We consider a system 
\begin{equation}
   \dot{x}=\begin{bmatrix}
    0 & 1\\-0.64 & -0.16
    \end{bmatrix}x+\begin{bmatrix}
    0\\1
    \end{bmatrix}u. 
\end{equation} 
$D=20I$ and $E=1$ are chosen resulting in an overdamped closed-loop system. The states were initialised along the first eigen mode. (i) Using the soft-constrained method, the recovered weights were far away from the true weights. The rank of $Q_p$ was found to be 3 meaning it was not full rank. Hence, the IOC problem is not solvable using the soft-constrained method for these trajectories. (ii) Using the hard-constrained method, the recovered weights were far away from the true weights. $Y_a$ was not full rank, leading to $W_{\psi}$ be of rank 1, that is, less than $k-1$, where $k=3$, in this case.

\textbf{Example 3}: We consider a system 
\begin{equation}
   \dot{x}=\begin{bmatrix}
    -0.5 & 1 \\ -1.895588534068740 & 1.393261306481120
\end{bmatrix}x+\begin{bmatrix}
    0 \\1
\end{bmatrix}u 
\end{equation} 
$D=0.0019I$ and $E=1$ are chosen resulting in a closed-loop under-damped system. The initial states were chosen as $x_0=\begin{bmatrix}
    1 & 0.091542807321846
\end{bmatrix}^T$ and $\text{Re}\{\hat{\lambda}+\sigma \}<0$.  (i) Using the soft-constrained method, the recovered weights were close to the true weights because $Q_p$ was full rank. (ii) Using the hard-constrained method, the recovered weights were far from the true weights because $\delta_1\mu_2\delta_2\mu_1< 0$ for this example and the initial conditions result in $H_rN$ and $H_cN$ in (\ref{hrn}) and (\ref{hcn}) to be linearly dependent.

\section{Conclusion}

The paper analyses the solvability of the IOC problem for two existing minimum principle based IOC methods. For the soft-constrained method, the observability matrix of the time-varying secondary problem was analysed while for the hard-constrained method, the rank of a matrix associated with the secondary optimization problem was analysed. Using trajectories obtained from over-damped and critically damped primary systems initialised at a single real mode, it was found that the IOC problem does not converge to the true solution for both methods. For trajectories obtained from under-damped closed-loop primary systems, the IOC problem converges to the true solution using the soft-constrained method. For the hard-constrained method, using under-damped closed-loop trajectories, the true weights can be recovered for any non-zero initial conditions if (i) the real part of sum of the closed-loop and open-loop eigenvalues is negative and (ii) the open-loop and closed-loop system parameters have a certain structure. The mathematical analysis was supported by appropriate presentation of specific numerical examples. The results give an understanding as to what kinds of trajectories depending on the initial conditions of the closed-loop system and system dynamics, of the primary optimal control problem, will result in either the success or the failure of the soft and the hard-constrained IOC methods based on the minimum principle. Current analytical and numerical results seem to suggest that for under-damped closed loop systems, the soft-constrained approach allows solvability of the IOC problem  whereas with the hard-constrained approach, the IOC problem is not solvable under certain conditions.






\section*{References}

\begin{thebibliography}{00}


\bibitem{c20} M. Johnson, N. Aghasadeghi, and T. Bretl, Inverse optimal control for deterministic continuous-time nonlinear systems, in $52^{nd}$ IEEE Conference on Decision and Control, IEEE, 2013, pp. 2906–2913.

\bibitem{c10a} T. L. Molloy, J. Inga, M. Farad, J. J. Ford, T. Perez and S. Hohmann, “Inverse Open-Loop Noncooperative Differential Games and Inverse Optimal Control,” IEEE Transactions on Automatic Control, vol. 65, no. 2, pp. 897-904, 2020. 

\bibitem{ca} Naidu, D. Subbaram, "Optimal control systems", CRC press, 2002.

\bibitem{c1} R. E. Kalman, “When is a linear control system optimal?” Journal of Basic Engineering, vol. 86(1), pp. 51–60, 19

\bibitem{c2} P. Abbeel and A. Y. Ng, “Apprenticeship learning via inverse reinforcement learning,” in Proceedings of the twenty-first international conference on Machine learning, 2004, pp. 1

\bibitem{c3} A.Y.Ng and S. J. Russell, “Algorithms for inverse reinforcement learning,” in International Conference on Machine Learning, vol. 1, 2000, p. 2

\bibitem{c4} N. Ratliff, B. Ziebart, K. Peterson, J. A. Bagnell, M. Hebert, A. K. Dey, and S. Srinivasa, “Inverse optimal heuristic control for imitation learning,” in Artificial Intelligence and Statistics, PMLR, 2009, pp. 424–43

\bibitem{c5}  Finn, S. Levine, and P. Abbeel, “Guided cost learning: Deep inverse optimal control via policy optimization,” in International conference on machine learning, PMLR, 2016, pp. 49–58.

 \bibitem{c19} K. Mombaur, A. Truong, and J. P. Laumond. "From human to humanoid locomotion—an inverse optimal control approach." Autonomous robots 28, no. 3 (2010): 369-383.

\bibitem{c7a} A. Keshavarz, Y. Wang and S. Boyd,“Imputing a convex objective function,” in 2011 IEEE international symposium on intelligent control, pp. 613-619, 2011.

\bibitem{c9} T. L. Molloy, J. J. Ford, and T. Perez, “Inverse noncooperative dynamic games,” IFAC PapersOnLine, vol. 50, no. 1, pp. 11788–11 793, 2017.

\bibitem{c10} T. L. Molloy, J. J. Ford, and T. Perez, “Finite-horizon inverse optimal control for discrete-time nonlinear systems,” Automatica, vol. 87, pp. 442–4, 2020.




\bibitem{c10b} W. M. Wonham, “On a Matrix Riccati Equation of Stochastic Control,”  SIAM Journal on Control, vol. 6, no. 4, pp. 681-697, 1968.



\bibitem{c12} R. Self, M. Abudia, S. N. Mahmud and R. Kamalapurkar, "Model-based inverse reinforcement learning for deterministic systems. Automatica", 140, 110242, 2022.


\bibitem{c14} R. Self, K. Coleman, H. Bai and R. Kamalapurkar, "Online observer-based inverse reinforcement learning", IEEE Control Systems Letters, 5(6), pp.1922-1927, 2020.

\bibitem{c15} B. Lian, W. Xue, F. L. Lewis, and T. Chai, "Online inverse reinforcement learning for nonlinear systems with adversarial attacks", International Journal of Robust and Nonlinear Control, 31(14), 6646-6667, 2021.


\bibitem{c18} B. Lian, W. Xue, F.L. Lewis, and T. Chai, "Inverse reinforcement learning for multi-player noncooperative apprentice games. Automatica", 145, p.110524, 2022.

\bibitem {c22} R. S. Bucy, "Global Theory of the Riccati Equation", Journal of Computer and System Sciences, 1, 349-361, 1967.

\bibitem{c23} H. d’Angelo, "Linear time-varying systems: analysis and synthesis", Allyn and Bacon, 1970.

\bibitem{c24} F.L. Lewis, D. Vrabie, and V. L. Syrmos," Optimal control", John Wiley \& Sons, 2012.

\bibitem{c25} T. Kailath, "Linear systems". Vol. 156. Englewood Cliffs, NJ: Prentice-Hall, 1980.

\end{thebibliography}
\end{document}